\documentclass[12pt]{amsart}
\usepackage[english]{babel}
\usepackage{amsmath}
\usepackage{amsthm}
\usepackage{amssymb}
\usepackage{mathrsfs}
\usepackage{enumerate}
\usepackage[notcite, final, notref]{showkeys}
\usepackage{dsfont}
\def\eop{\hfill\square}
\usepackage[dvips]{color}
\newtheorem{theorem}{Theorem}[section]

\newtheorem{lemma}[theorem]{Lemma}

\theoremstyle{definition}
\newtheorem{remark}[theorem]{Remark}

\def\E{\bf E}

\usepackage{xcolor}

\title[]{A Sparse Representation of Random Signals}

\author[T. Qian]{Tao Qian*}
\address{Tao QIAN, Macao Center for Mathematical Sciences \\
Macau University of Science and Technology\\
Macau}
\email{tqian@must.edu.mo}
\thanks{*Corresponding author.\\
Funded by The Science and Technology Development Fund, Macau SAR (File no. 0123/2018/A3)}

\def\C{\bf C}

\begin{document}
\maketitle
\begin{abstract} Studies of sparse representation of deterministic signals have been well developed. Amongst there exists one called adaptive Fourier decomposition (AFD) established through adaptive selections of the parameters defining a Takenaka-Malmquist system in one-complex variable. The AFD type algorithms give rise to sparse representations of signals of finite energy. The multivariate generalization of AFD is one called pre-orthogonal AFD (POAFD), the latter being established with the context Hilbert space possessing a dictionary. The purpose of the present study is to generalize both AFD and POAFD to random signals. We work on two types of random signals. One is those expressible as the sum of a deterministic signal with an error term such as a white noise; and the other is, in general, as mixture of several classes of random signals obeying certain distributive law. In the first part of the paper we develop an AFD type sparse representation for one-dimensional random signals by making use analysis of one complex variable. In the second part, without complex analysis, we treat multivariate random signals in the context of stochastic Hilbert space with a dictionary. Like in the deterministic signal case the established random sparse representations are powerful tools in practical signal analysis.
\end{abstract}

\bigskip


\def\D{\bf D}
\def\C{\bf C}
\bigskip

\date{today}
\tableofcontents
\def\R{\bf R}
\section{Introduction }
If $F$ is a complex-valued signal in $[0,2\pi)$ with finite energy, then it
can be expanded into its $L^2([0,2\pi))$-convergent Fourier series:
\[ F(t)=\sum_{k=-\infty}^\infty c_ke^{ikt}.\]
To make convenient use of complex analysis we alter the notation and denote it as $f(e^{it})=F(t).$
Then the Plancherel Theorem asserts the relation $\|f\|^2=\sum_{-\infty}^\infty |c_k|^2,$ where the $L^2$-norm is one with respect to the inner product
 \[ \langle f,g\rangle=\frac{1}{2\pi}\int_0^{2\pi}f(e^{it})\overline{g}(e^{it})dt.\]
 The Plancherel relation infers that $c_k$ tends to zero and therefore the complex-valued functions
\[ f^+(z)=\sum_{k=0}^\infty c_kz^k \quad {\rm and}\quad f^-(z)=\sum_{k=-1}^{-\infty} c_kz^k\]
are analytic in  $\D$ and in $\C\setminus \overline{\D},$ respectively, where $\D$ stands for the open unit disc in the complex plane $\C.$
Restricted to the unit circle, in the $L^2$-convergence sense, we define
\[ f^+(e^{it})\triangleq\sum_{k=0}^\infty c_ke^{ikt}\]
as the \emph{analytic signal} associated with $f.$ Denote by $H$ the Hilbert transform operator on the circle defined by
\[ Hf(e^{it})=\sum_{k=-\infty}^\infty (-i){\rm sgn}(k)c_ke^{ikt},\]
where ${\rm sgn}(k)=k/|k|$ when $k\ne 0;$ and ${\rm sgn}(0)=0.$ Restricted on the circle we have
$f^\pm=\frac{1}{2}(f+iHf\pm c_0).$
The non-tangential boundary limit of $f^+(z)$ as $z\to e^{it}$ coincides with
the above defined $L^2$-limit $f^+(e^{it}).$
 To be practical we assume that the test functions $f$ are real-valued. Then $c_{-n}=\overline{c}_n,$ and, as a consequence, \[ f(e^{it})=2{\rm Re}\{f^+(e^{it})\}-c_0.\]
Due to the above relation, harmonic analysis of
signals $f$ of finite energy can be reduced
to complex analysis of the associated analytic signals $f^+.$ Since $f^+$ is the boundary limit of the analytic function $f^+(z)$ in $\D,$  complex analytic methods are available for $f^+.$  The totality of all such analytic functions $f^+(z)$ in the disc with finite energy constitute the function space
\begin{eqnarray}\label{alter} H^2({\D})&\triangleq&\{f:{\D}\to {\C}\ |\ f\ {\rm is\ analytic\ and}\ f(z)=\sum_{k=0}^\infty c_kz^k\ {\rm with}\ \sum_{k=0}^\infty |c_k|^2<\infty\}\nonumber \\
&=&\{f:{\D}\to {\C}\ |\ f\ {\rm is\ analytic\ and}\ \sup_{0<r<1}\int_0^{2\pi}|f(re^{it})|^2dt<\infty\},\end{eqnarray}
called the (\emph{complex}) \emph{Hardy $H^2$-space} in the unit disc. There exist other complex Hardy spaces having more or less parallel theories: The Hardy space idea to study functions may be extended to signals defined on the whole real line $\R,$ to those defined on manifolds in the higher dimensional complex spaces ${\C}^d$ with the several complex variables setting (e.g., the Hardy spaces on tubes \cite{SW}), or to those in the real-Euclidean spaces ${\R}^d$ in the Clifford algebra setting (the conjugate harmonic systems, \cite{SW}. And the functions under study can have scalar-, or complex-, or vector-, and even matrix-values (\cite{ACQS1,ACQS2}), etc. All obey the same philosophy. We will only take the context $H^2({\D})$ as an example to explain the adaptive Fourier decomposition (AFD) theory. In below we
 often abbreviate $H^2({\D})$ as $H^2.$ The Hardy space $H^2({\D})$ has several equivalent characterizations that are not of interest in this paper. The disc case corresponds to signals defined in a compact interval on the line. That is the mathematical formulation for periodic signals. In the first half of this paper we mainly concentrate in a stochastic-lization of the Hardy space. In the Hardy space the adaptive Fourier decomposition, abbreviated as AFD or Core-AFD, was first established in \cite{QWa}. We note that AFD on the disc and on the real line heavily depends on the concepts Blaschke product and Takenaka-Malmquist (TM) system. The related delicate analysis achieves AFD as a beautiful piece of mathematical work well fitting into the Beurling-Lax structure of the classical complex Hardy spaces.
 In many analytic function spaces, functions with the same role as Blaschke products are not available. Pre-orthogonal AFD (POAFD) provides a replacement of AFD in Hilbert spaces that usually do not have Blaschke product-like functions or TM systems, nor explicit orthogonal function systems. The Hilbert space setting are, in particular, for multivariate signals. We leave the POAFD method to be studied in the second half of this paper in which we formulate stochastic POAFD in the setting of stochastic Hilbert space with a dictionary.\\
\def\N{\mathcal{N}}
In practice one encounters random signals: Signals are mostly corrupted with noise or together with measurement errors, or, as an alternative type, belong to a collection consisting of signals in different classes obeying certain distribution law. A practical formulation then should be a real-valued function $F(t,w),$ where for a fixed probabilistic sample point $w\in \Omega$ the function $F(\cdot,w)$ is a deterministic signal; meanwhile for each point $t$ in the time domain or the space domain the function $F(t,\cdot)$ is a random variable. We call such signals \emph{random signals} (RSs). To formulate the corresponding stochastic Hardy space theory in the case $t\in [0,2\pi)$ we denote $F(t,w)=f(e^{it},w),$ and write the
trigonometric expansion of the latter as
\[ f(e^{it},w)=\sum_{k=-\infty}^\infty c_k(w)e^{ikt}=[\sum_{k=-\infty}^\infty c_k(w)z^k]_{z=e^{it}},\  {\rm where} \ \ c_k(w)=\frac{1}{2\pi}\int_0^{2\pi}f(e^{is},w)e^{-iks}ds. \]
The Plancherel Theorem gives
\[  \|f(\cdot,w)\|^2=\sum_{k=-\infty}^\infty |c_k(w)|^2.\]
In our study it is reasonable to impose the condition
\begin{eqnarray}\label{bdd}
\left[E_w\|f(\cdot,w)\|^2\right]^{\frac{1}{2}}=
\left(\sum_{k=-\infty}^\infty E_w|c_k(w)|^2\right)^{\frac{1}{2}}<\infty,
\end{eqnarray}where $E_w$ stands for the mathematical expectation in the underlying probability space. In the whole paper the underlying probability space, $(\Omega,\mu), w\in \Omega,$ is not specified. The theory to be developed will be valid for any but fixed probability space. The quantity in (\ref{bdd}) is called the \emph{energy expectation norm} (EE-Norm) of $f,$ denoted as $\|f\|_{\mathcal{N}}.$ Set,
 \begin{eqnarray}\label{normal}
 L^2_w({\partial \D},\Omega)=\{f:\partial {\D}\times\Omega\to {\C}\ |\ f \ {\rm is\ a\ real{\mbox -}valued\ RS,\ and}\ \|f\|_{\mathcal{N}}<\infty\},
 \end{eqnarray}

called \emph{the space of} \emph{random signals of} \emph{finite energy}. $L^2_w({\partial \D},\Omega)$ is written briefly as $\mathcal{N}.$ The RSs in $L^2_w(\partial {\D},\Omega)$ are called \emph{normal random signals}, or normal RSs.  The space $\mathcal{N}$ is a Hilbert space under the inner product induced from the EE-norm. A normal RS is almost surely a signal of finite energy in $t.$
In below we will keep the inner product notation $\langle \cdot,\cdot\rangle$ only for the inner product of the time-domain-space $L^2(\partial {\D}).$ We note that although we have the restriction that $f\in \N$ is real-valued, for many results that we deduce below the restriction is unnecessary. \\

Similarly to the deterministic case we will concentrate in studying $\lq\lq$a half"
 of the space $\mathcal{N},$ consisting of the RSs with expansions in the spectrum range $k=0,1,\cdots,$
\[ f^+(e^{it},w)=\sum_{k=0}^\infty c_k(w)e^{ikt}, \quad {\rm satisfying }\quad \sum_{k=0}^\infty
E_w(|c_k(w)|^2)<\infty.\]
As a consequence, almost surely,
\[ \sum_{k=0}^\infty |c_k(w)|^2<\infty.\]
The last condition implies that $|c_k(w)|$ is a bounded sequence. Hence, almost surely,
\[ f^+(z,w)=\sum_{k=0}^\infty c_k(w)z^k\] is an analytic function in ${\D}.$ The boundary limits exist in the a.e. pointwise, and in the $L^2$-convergence sense as $r=|z|\to 1,$  and
\[ f(e^{it},w)=2{\rm Re}\{f^+(e^{it},w)\}-c_0(w).\]
On the boundary $\partial \D$ the projection $f^+,$ apart from being obtained through the Taylor series expansion, can also be obtained through the (circular) Hilbert transform $H:$
\begin{eqnarray}\label{given by} f^+(e^{it},w)=\frac{1}{2}(f(e^{it},w)+iHf(e^{it},w)+c_0),\end{eqnarray}
where for any $f(e^{it},w)=\sum_{k=-\infty}^\infty c_k(w)e^{ikt},$ in terms of the signum function ${\rm sgn}{(k)}=k/|k|, k\ne 0$ and ${\rm sgn}{(0)}=0,$ the Hilbert transform $Hf$ can be defined through the Fourier multiplier form or, alternatively, through the singular integral form,
\begin{eqnarray*}\label{Hilbert}
 Hf(e^{it},w)&\triangleq&\sum_{k=-\infty}^\infty (-i){\rm sgn}(k)c_k(w)e^{ikt}\\
 &=&\frac{1}{\pi}{\rm v.p.}\int_{-\infty}^\infty\cot\left(\frac{s}{2}\right)f(e^{i(t-s)},w)ds.
 \end{eqnarray*}
 It is obvious that Hilbert transformation maps real-valued functions to real-valued functions.

By using Hilbert transform, the study of normal RSs on the circle is reduced to the study of the half series constituting the Hardy space in the disc. The philosophy and related techniques may be found in the literature  \cite{QWa,Qian2010,Q2D,CQT} and the references therein.

We define the \emph{stochastic Hardy space} as follows (with the superscript $\lq\lq$${}^+$" suppressed), denoted
\def\N{\mathcal{N}}
\begin{eqnarray}
H^2_w({\D})&=&\{f:{{\D}}\times \Omega\to {\C}\ |\ f(z,w)\ {\rm is\ a.s.\  analytic\ in}\ z \ {\rm and}\nonumber \\
& &\qquad f(z,w)=\sum_{k=0}^\infty c_k(w)z^k\ {\rm with}\ \|f\|_{\N}^2=\sum_{k=0}^\infty E_w|c_k(w)|^2<\infty\}.
\end{eqnarray}
Denoted by $H^2_w(\partial \D)$ the function space consisting of the non-tangential boundary limits of the RSs in $H^2_w({\D}).$ On the boundary $\partial \D$ the space $H^2_w(\partial \D)$ is a proper closed subspace of $\mathcal{N},$ being easily seen in view of the Fourier series expansions of the functions. \\

\begin{remark}
There is an analogous formulation for random signals in the whole time range, where the corresponding space $\N=L^2_w({\R},\Omega)$ is defined as
\begin{eqnarray*}\label{normal}L^2_w({\R},\Omega)=\{f:{\R}\times\Omega\to {\C}\ |\ f \ {\rm is\ a\ real{\mbox -}valued\ RS,\ and}\ \|f\|^2_{\mathcal{N}}=E_\omega\|f_\omega\|^2_{L^2({\R})}<\infty\},\end{eqnarray*}
and the random Hardy space is defined as
\begin{eqnarray*}
H^2_w({\C}^+)&=&\{f:{{\C}^+}\times \Omega\to {\C}\ |\ f(z,w)\ {\rm is\ a.s.\  analytic\ in}\ z \ {\rm and}\nonumber \\
& &\quad f(z,w)=\int_0^\infty F(\xi,w)e^{-i\xi z}d\xi\ {\rm with}\ \|f\|_{\N}^2=\int_0^\infty E_w|F(\xi,w)|^2d\xi<\infty\}.
\end{eqnarray*}
This formulation is related to Laplace transform. It also have a well developed Hilbert transform theory.
\end{remark}

The purpose of this study is to develop stochastic adaptive Fourier decompositions (SAFDs) for analyzing random signals of two types: noised signals and classes of signals. We will accordingly develop two models of stochastic AFD (SAFD), namely SAFDI and SAFDII. What makes the complex analysis method a great power in the single variable case is the Cauchy theory, involving the Cauchy theorem and the Cauchy formula, the latter reproduces the function values of an analytic function inside its domain by using its data on the boundary of the domain. Thanks to Cauchy's theory and availability of Blaschke products, AFD and SAFD stand as beautiful mathematical theory with profound relations to classical complex analysis, and in particular to rational function approximation. A natural generalization of the analytic function theory to multivariate functions would be one for reproducing kernel Hilbert spaces. A close examination shows that although it offers convenience, a reproducing kernel is, in fact, unnecessary, but a dense and complete set, for further developing the theory. A Hilbert space with a dictionary is a more general concept than a reproducing kernel Hilbert space. In the later half of this paper we extend the sparse approximation theory for analytic RSs in one-complex variable to \emph{stochastic Hilbert space} with a dictionary. The necessity of developing a theory in the general Hilbert space context rests in the tendency of studying various multivariate random signals,
 in which there may not exist analyticity properties as used in the classical Hardy spaces case.\\

The writing plan is as follows. In \S 2 with the stochastic Hardy space setting we establish two types of sparse approximations, SAFDI and SAFDII, for treating two categories of  analytic RSs: One is for noised deterministic signals, and the other is a collection of several classes of signals obeying certain probability distribution law. In \S 3 we extend the theory to the context of stochastic Hilbert space with a dictionary treating also the same two categories of RSs, and develop their respective sparse approximations, SPOAFDI and SPOAFDII.

For the reader's convenience we give the following abbreviations list:\\

AFD: adaptive Fourier decomposition (for deterministic signals in the classical Hardy spaces consisting of analytic signals of finite energy on the boundary, associated with a Blaschke product structure)\\

BVC: boundary vanishing condition\\

MSP: maximal selection principle\\

POAFD: pre-orthogonal adaptive Fourier decomposition (Applicable for Hilbert spaces with a dictionary satisfying BVC)\\

SBVC: stochastic boundary vanishing condition\\

RS: random signal\\

Normal RS: normal random signal, or a signal in the space (\ref{normal})\\

$\mathcal{N}$: the Hilbert space consisting of normal RSs\\

$H^2_w({\D}):$ the stochastic Hardy space on the disc, corresponding to $c_k(w)=0$ for $k<0$\\

$H^2_w(\partial {\D}):$ the space of the functions as boundary limits of those in $H^2_w({\D})$ defined on $\partial {\D}$\\

SHS: a stochastic Hilbert space, or a Hilbert space of RSs possessing finite variation\\

SAFD, SAFDI, SAFDII: stochastic AFDs (SAFDs) are divided into two types: the type I, SAFDI, is for the RSs that are expressible as a deterministic signal corrupted with a noise of small $\mathcal{N}$-norm; the type II, SAFDII, is for a general stochastic Hardy space.\\

SPOAFD, SPOAFDI, SPOAFDII: stochastic POAFDs (SPOAFDs) in SHS consist of two types; the type I, SPOAFDI, is for the RSs being expressible as noised signals; the  type II, SPOAFDII, is for any general SHS.

\section{Stochastic AFDs}

In the deterministic signal analysis AFD is a sparse approximation methodology using a suitably adapted Takenaka-Malmquist (TM) system. In the classical Hardy space formulation it well fits with the Beurling-Lax Theorem, where any specific function belongs to a backward-shift-invariant subspace in which the function is the limit of a fast converging TM series. The AFD type expansions have found
many applications in signal and image analysis as well as in system identification (see, for instance, \cite{LZQ,WQLG,Mi1,Mi2}). With the stochastic Hardy space defined in \S 1 we present two types of AFD-like expansions, called \emph{stochastic AFDs} (SAFDs), being used for different purposes in application. Before studying SAFDs we develop some aspects in relation to the Hardy space projections of the normal RSs.

\subsection{Properties of Hardy Space Projections of Random Signals}
\def\r{\tilde{r}}
Normal RSs $f(e^{it},w)$ can all be represented into the form
\begin{eqnarray}\label{rTypeI} f(e^{it},w)=\tilde{f}(t)+\tilde{r}(e^{it},w),\end{eqnarray}
where $\tilde{f}=E_wf.$ The difference $\r$ is
sometimes called the \emph{remainder RS}.
In this section we reduce the analysis of ordinary normal RSs to that of the analytic normal RSs. The philosophical support of this methodology is the relation (\ref{given by}). Given by the following two theorems, the Hardy space projections $f^+, \tilde{f}^+$ and $\r^+$ possess the most favorable properties own by those from which they are projected.

 \begin{theorem}\label{basic}
 If $f\in \N,$ then $\tilde{f}\in L^2(\partial {\D}), \r\in \N, E\r=0.$  In writing
 \[f(e^{it},w)=\sum_{k=-\infty}^\infty c_k(w)e^{ikt}\quad {\rm and}\quad \r(e^{it},w)=\sum_{k=-\infty}^\infty d_k(w)e^{ikt},\]
 there hold
 \[ \tilde{f}(e^{it})=\sum_{k=-\infty}^\infty (E_wc_k)e^{ikt},\]
 where $E_wc_k=E_w(c_k(w)),$ and,
 \[ d_k(w)=c_k(w)-E_wc_k, \quad E_wd_k=0, \quad k=0,\pm 1, \pm 2\cdots\]
 The Hardy space projections $f^+, \tilde{f}^+, \r^+,$ respectively, belong to $H^2_w(\partial {\D}), H^2(\partial \D),$ and $H^2_w(\partial \D).$ There hold
\[ \{E_wf\}^+=E_w\{f^+\} \quad {\rm and}\quad \|{\r}^+\|_{\N}=\frac{\|\r+d_0\|_{\N}}{\sqrt{2}}.\]
 \end{theorem}
 \noindent{\bf Proof} We note that
\begin{eqnarray}\label{follow}
\left(\sum_{k=-\infty}^\infty|E_w(c_k(w))|^2\right)^{1/2}&\leq&
E_w\left[\left(\sum_{k=-\infty}^\infty|c_k(w)|^2\right)^{1/2}\right]\quad ({\rm Minkovski's \ inequality})\nonumber \\
&\leq& \left[E_w(\sum_{k=-\infty}^\infty|c_k(w)|^2)\right]^{1/2}[E_w(1)]^{1/2}\quad
({\rm H\ddot{o}lder's \ inequality})\nonumber \\
&=&\left[\sum_{k=-\infty}^\infty E_w(|c_k(w)|^2)\right]^{1/2}[E_w(1)]^{1/2}\nonumber \\
&=& \|f\|_{\N} <\infty.
\end{eqnarray}
Then the Riesz-Fisher Theorem asserts that
\[g(e^{it})=\sum_{k=-\infty}^\infty E_w(c_k(w))e^{ikt}\in L^2(\partial \D).\]
Now we show $\tilde{f}=g.$ Denote $f_n(e^{it},w)=\sum_{|k|\leq n}c_k(w)e^{ikt}.$ Then
$E_wf_n(e^{it},w)=\sum_{|k|\leq n}E_w(c_k)e^{ikt}.$ Similarly to the reasonings in proving (\ref{follow}), there follow
\begin{eqnarray*}
\|E_wf-E_wf_n\|&=& \|E_w(f-f_n)\|\\
&\leq& E_w\|f-f_n\|\\
&\leq& \left(E_w\|f-f_n\|^2\right)^{1/2}\\
&=&\|f-f_n\|_{\N}\\
&=&\left(\sum_{|k|>n}E_w(|c_k(w)|^2)\right)^{1/2}\\
&\to& 0, \quad {\rm as}\ n\to\infty.
\end{eqnarray*}
Since the linear functional of the $m$-th Fourier coefficient, $C_m,$ is continuous, there follows
\[ C_m(E_wf)=\lim_{n\to \infty}C_m(E_wf_n)=E_w(c_m).\]
This shows that $E_wf=g\in L^2(\partial \D)$ and is with the Fourier expansion
\[ \tilde{f}=\sum_{k=-\infty}^\infty E_w(c_k(w))e^{ikt}\in L^2(\partial \D).\]
It then follows
\begin{eqnarray}\label{zero mean} E_w(\r(e^{it},w))=E_wd_k=0,\quad \forall t\in [0,2\pi) \ {\rm and}\ k=0,\pm 1\cdots\end{eqnarray}
As a consequence of (\ref{zero mean}), we have the orthogonality
\begin{eqnarray}\label{conse} E_w(|\tilde{f}(e^{it})+\r(e^{it},w)|^2)=|\tilde{f}(e^{it})|^2+E_w(|\r(e^{it},w)|^2).\end{eqnarray}
Thus, by taking the integration with respect to $dt,$ we have the finiteness of the $\N$-norm of  $\r:$
\begin{eqnarray}\label{bddE}
\|{\r}\|_{\N}^2=\|f\|_{\N}^2-\|\tilde{f}\|_{L^2(\partial \D)}^2<\infty. \end{eqnarray}
 To compute the $\N$-norm of $\tilde{r}^+,$ by taking into account
 $d_k=\overline{d}_{-k},$
we have
\begin{eqnarray*}
\|r^+\|_{\N}^2=E_w\int_0^{2\pi}|r^+(e^{it},w)|^2dt=\sum_{k=0}^\infty
E_w|{d}_k(w)|^2=\frac{\|\tilde{r}+{d}_0\|_{\N}^2}{2}.
\end{eqnarray*}
 The proof of the theorem is complete.$\eop$\\

 A particular example is that the remainder $\r$ in question is the white Gaussian noise $N(0,\sigma^2),$ when the relation (\ref{bddE}) becomes
 \[ E_w(|{\r}(e^{it},w)|^2)= \sigma^2, \quad \forall t\in [0,2\pi).\]
\def\N{\mathcal{N}}

 We would be interested in properties imposed to the remainder RS $\r$ not as special as white noise. What having in mind are weakly stationary, and further, ergodic RS $\r.$ Since we already have $E_w\r=0,$ recall that if
 the autocorrelation function of $\r$ (coinciding with the autocovariance function of $f$ itself) depends only on the time difference, that is, if there holds for some deterministic signal $\tilde{r}_1,$
 \begin{eqnarray}\label{stationary}
 \tilde{\gamma} (t,s)=E_w({\r}(e^{it},w){\r}(e^{is},w))\triangleq
 \tilde{r}_1(s-t),\end{eqnarray}
 then $\tilde{r}$ is called a \emph{weakly stationary RS}.

Recall that a weakly stationary RS, say $x(t,w),$ is \emph{weakly ergodic} if and only if
\begin{eqnarray}\label{er1}
E_wx=\lim_{T\to \infty}\frac{1}{2T}\int_{-T}^T x(t,w)dt, \quad {\rm a.s.},
\end{eqnarray}
and
\begin{eqnarray}\label{er2}
E_w(x(t,w)\overline{x}(t-\tau,w))=\lim_{T\to \infty}\frac{1}{2T}\int_{-T}^T x(t,w)\overline{x}(t-\tau,w)dt, \quad {\rm a.s.}
\end{eqnarray}

The condition (\ref{er1}) implies that the common quantity of the LHS and the RHS of the equality (\ref{er1}) is almost surely a constant. The condition (\ref{er2}) implies that the common quantity of the LHS and the RHS of the equality (\ref{er2}) is almost surely a function of the time difference $\tau.$

Since $E_w\r=0,$ the relation (\ref{rTypeI}) implies that, under the condition $\r$ being weakly stationary, $f$ is weakly stationary if and only if $\tilde{f}$ is almost surely a constant function; and, $f$ is weakly ergodic if and only if $\tilde{f}$ is almost surely the zero function.
A pure random variable $f(w)$ is stationary. If it is further ergodic, then it has to be a constant almost surely. This observation together with the above one hint that it would be necessary to assume $d_0=0$ when discussing stationarity and ergodicity of RSs.

\begin{theorem}\label{projection} Under the assumptions of Theorem \ref{basic}, if further $d_0=0,$ a.s., then weak stationarity of $\r$ implies weak stationarity of $\r^+;$  and, weak ergodicity of $\r$ implies weak ergodicity of $\r^+.$
\end{theorem}
We need first prove the following lemma.

\begin{lemma}\label{la1} The Hilbert transform
$H$ and the expectation operator $E_w$ are commutable.
 \end{lemma}
 \noindent{\bf Proof}
 As proved in the beginning of the proof of Theorem \ref{basic}, the series
\[ \sum_{k=-\infty}^\infty |E_wc_k|^2\]
is convergent. It implies that
\[\sum_{k=-\infty}^\infty (-i){\rm sgn}(k)(E_wc_k)e^{ikt}\]
is  convergent in $L^2(\partial {\D}).$ This implies, as in the proof of Theorem \ref{basic},
\[ E_w\sum_{k=-\infty}^\infty (-i){\rm sgn}(k)c_k(w)e^{ikt}=\sum_{k=-\infty}^\infty (-i){\rm sgn}(k)(E_wc_k)e^{ikt}.\] Hence,
\begin{eqnarray*} (E_wH)f(e^{it})&=&E_w\sum_{k=-\infty}^\infty (-i){\rm sgn}(k)c_k(w)e^{ikt}\\
&=&\sum_{k=-\infty}^\infty (-i){\rm sgn}(k)E_w(c_k(w))e^{ikt}\\
&=& H(E_wf)(e^{it}).
\end{eqnarray*}
The proof is complete.$\eop$\\

The limits of function sequences (in variable $w$) on the RHDs of (\ref{er1}) and (\ref{er2}) are understood as pointwise. In the proof of the following theorem the convergence is assumed to be in suitable function spaces on which the expectation functional and the Hilbert transformation operator are continuous.\\

\noindent {\bf Proof of Theorem \ref{projection}}
We first show that if the autocorrelation function of ${\r}$ is a function of, merely, the time difference, then that of $\r^+$ is the same.  For this goal we note the fact that weakly stationarity of $\r$ implies the orthogonality under the expectation operation: $E_w(\tilde{d}_k\overline{\tilde{d}}_l)=\delta_k(l),$ where $\delta_k(l)$ is the Dirac Delta function. Expanding ${\r}^+(e^{it},w)$ and ${\r}^+(e^{is},w)$ into their respective Fourier series, the noted fact implies that
\[ E_w({\r}^+(e^{it},w)\overline{{\r}^+}(e^{is},w))=\sum_{k=0}^\infty E_w(|c_k(w)|^2)e^{ik(t-s)},\]
depending only on the time difference.
Next, we assume that $\r$ is weakly stationary and weakly ergodic. We first show that the expectation is ergodic. By invoking the commutativity between $H$ and $E_w$ proved in Lemma \ref{la1}, and the property that the Hilbert transform $H$ annihilates constant functions, we have, almost surely,
\begin{eqnarray*}
E_w(\r^+)&=&\frac{1}{2}E_w(\r+iH\r+d_0)\\
&=&\frac{1}{2}(E_w\r+iE_wH\r)\\
&=&\frac{1}{2}E_w\r+i\frac{1}{2}HE_w\r\\
&=&\frac{1}{2}E_w\r.
\end{eqnarray*}
On the other hand,
\begin{eqnarray*}
\lim_{T\to \infty}\frac{1}{2T}\int_{-T}^T\r^+(e^{it},w)dt&=&\frac{1}{2}\lim_{T\to \infty}\frac{1}{2T}\int_{-T}^T\r(e^{it},w)dt+\frac{i}{2}\lim_{T\to \infty}\frac{1}{2T}\int_{-T}^TH\r(e^{it},w)dt\\
&=&\frac{1}{2}\lim_{T\to \infty}\frac{1}{2T}\int_{-T}^T\r(e^{it},w)dt+\frac{i}{2}\lim_{T\to \infty}H(\frac{1}{2T}\int_{-T}^T\r(e^{it},w)dt)\\
&=&\frac{1}{2}E_w\r+\frac{i}{2}H(E_w\r)\\
&=&\frac{1}{2}E_w\r.
\end{eqnarray*}
Therefore the expectation is ergodic. To show that the autocorrelation is also ergodic we proceed similarly. We first write
\[
E_w(\r^+(e^{it},w)\overline{\r^+}(e^{i(t-s)},w))=
\frac{1}{4}E_w([\r(e^{it},w)+iH\r(e^{it},w)][\r(e^{i(t-s)},w)-iH\r(e^{i(t-s)},w)]).\]
The RHS of the last identity can be expressed as a complex linear combination of the following four terms:
\[E_w(\r(e^{it},w)\r(e^{i(t-s)},w)),\quad  E_w(H\r(e^{it},w)\r(e^{i(t-s)},w)),\]
\[E_w(\r(e^{it},w)H\r(e^{i(t-s)},w))\quad {\rm and} \quad E_w(H\r(e^{it},w)H\r(e^{i(t-s)},w)).\]
For the first term, due to the assumed ergodicity, we have
\[E_w(\r(e^{it},w)\r(e^{i(t-s)},w))=\lim_{T\to\infty}\frac{1}{2T}\int_{-T}^T
\r(e^{it},w)\r(e^{i(t-s)},w)dt.\]
We show that with each of the rest three terms the expectation operator may commute with the partial circular Hilbert transforms.  The commutativity then leads to the respective ergodicity. With a little abuse of the notation, temporarily denoting $$H_ug(e^{i(t-u)})=\frac{1}{\pi}{\rm v.p.}\int_{-\pi}^\pi \cot\frac{u}{2}g(e^{i(t-u)})du,$$ we have
\begin{eqnarray*}
E_w(H\r(e^{it},w)\r(e^{i(t-s)},w))&=&E_w(H_u(\r(e^{i(t-u)},w)\r(e^{i(t-s)},w)))\\
&=&H_u(E_w(\r(e^{i(t-u)},w)\r(e^{i(t-s)},w)))\\
&=&H_u(\lim_{T\to\infty}\frac{1}{2T}\int_{-T}^T\r(e^{i(t-u)},w)\r(e^{i(t-s)},w)dt)\\
&=&\lim_{T\to\infty}\frac{1}{2T}\int_{-T}^TH_u(\r(e^{i(t-u)},w)\r(e^{i(t-s)},w)dt)\\
&=&\lim_{T\to\infty}\frac{1}{2T}\int_{-T}^TH\r(e^{it},w)\r(e^{i(t-s)},w)dt.
\end{eqnarray*}
Similarly, we have
\[E_w(\r(e^{it},w)H\r(e^{i(t-s)},w))=\lim_{T\to\infty}\frac{1}{2T}\int_{-T}^T
\r(e^{it},w)H\r(e^{i(t-s)},w)dt.\]
For the last term we have
\begin{eqnarray*}
E_w(H\r(e^{it},w)H\r(e^{i(t-s)},w))&=&E_w(H_uH_v(\r(e^{i(t-u)},w)\r(e^{i(t-s-v)})))\\
&=&H_uH_v(E_w(\r(e^{i(t-u)},w)\r(e^{i(t-s-v)})))\\
&=&H_uH_v(\lim_{T\to\infty}\frac{1}{2T}\int_{-T}^T\r(e^{i(t-u)},w)\r(e^{i(t-s-v)})dt)\\
&=&\lim_{T\to\infty}\frac{1}{2T}\int_{-T}^TH_uH_v(\r(e^{i(t-u)},w)\r(e^{i(t-s-v)})dt)\\
&=&\lim_{T\to\infty}\frac{1}{2T}\int_{-T}^TH\r(e^{it},w)H\r(e^{i(t-s)},w)dt.
\end{eqnarray*}
The above exchange of limiting procedures may be first verified for a nice function class, and then extended to the whole underlying space through a density argument. A suitable linear combination of the four ergodic identities then gives rise to
 \[
E_w(\r^+(e^{it},w)\overline{\r^+}(e^{i(t-s)},w))=\lim_{T\to\infty}\frac{1}{2T}
\int_{-T}^T\r^+(e^{it},w)\overline{\r^+}(e^{i(t-s)},w)dt.\]
Ergodicity of the autocorrelation is thus proved. The proof of
Theorem \ref{projection} is complete.$\eop$\\

\subsection{The Type SAFDI: Taking Expectation First}

In this section we assume that $f(e^{it},w)$ is from $H^2_w(\D).$
Letting $\tilde{f}=E_w(f(e^{it},w)),$ we, as in the last section, have
\[ f(e^{it},w)=\tilde{f}(e^{it})+{\r}(e^{it},w).\]
The function $\tilde{f}$ is, in fact, in $H^2(\D).$ This is a consequence of Theorem \ref{basic}, or can be proved by the similar but the integral form inequalities as,
 for $r<1,$
\begin{eqnarray}\label{similar}
& &\left(\int_0^{2\pi}|E_wf(re^{it},w)|^2dt\right)^{1/2}\nonumber \\
&\leq& E_w\left[\left(\int_0^{2\pi}|f(re^{it},w)|^2dt\right)^{1/2}\right]\quad ({\rm Minkovski's\  inequality})\nonumber \\
&\leq& \left(E_w\int_0^{2\pi}|f(re^{it},w)|^2dt\right)^{1/2}E_w(1)^{1/2}\quad ({\rm Holder's \ inequality})\nonumber \\
&\leq& \|f\|_{\N}<\infty.\end{eqnarray}
We also note that, as a consequence of the last inequality, for a.s. $w\in \Omega,$ $f(re^{it},w)$ is a function in the classical complex Hardy space with the power series expansion
\[ f(re^{it},w)=\sum_0^\infty c_k(w)r^ke^{ikt}, \qquad r<1.\]
The type SAFDI is based on AFD of the deterministic signal $\tilde{f}.$ For the self-containing purpose we now go through a full AFD expansion of $\tilde{f}.$
We will be using the $L^2$-normalized Szeg\"o kernel on the circle:
\[ e_a(z)=\frac{\sqrt{1-|a|^2}}{1-\overline{a}z}, \quad a\in \D.\]
In $H^2({\D})$ it has the reproducing kernel property: For any $g\in H^2(\D),$
\[\langle g,e_a\rangle = \sqrt{1-|a|^2}g(a).\]
Let $f_1=\tilde{f}.$ For any $a\in {\D}$ we have the following identity as an orthogonal decomposition
\begin{eqnarray}\label{use1} \tilde{f}(z)=\langle f_1,e_a\rangle e_a(z) + f_2(z)\frac{z-a}{1-\overline{a}z},\end{eqnarray}
where $f_2$ is call the \emph{reduced remainder}, given by
\begin{eqnarray}\label{use2} f_2(z)=\frac{f_1(z)-\langle f_1,e_a\rangle e_a(z)}{\frac{z-a}{1-\overline{a}z}}\in H^2({\D}).\end{eqnarray}
Due to the orthogonalization we have
\begin{eqnarray}\label{use3}\|\tilde{f}\|_{H^2({\D})}=|\langle f_1,e_a\rangle|^2+\|f_2\|_{H^2({\D})}.\end{eqnarray}
Thus, the larger is the quantity $|\langle f_1,e_a\rangle|^2,$ the smaller is the energy of the reduced remainder $f_2.$
Although $\D$ is an open set it can be proved (see \cite{QWa}, for instance) that
\[ \sup\{|\langle f_1,e_a\rangle|^2\ |\ a\in {\D}\}\]
is attainable at a point of $\D.$ Hence, one practically selects
\[ a_1=\arg \max\{|\langle f_1,e_a\rangle|^2\ |\ a\in {\D}\}.\]
Such maximal selection is phrased as \emph{Maximal Selection Principle} (MSP) of the Hardy space (\cite{QWa}). The MSP is evidenced by the boundary vanishing condition (BVC) of the Szeg\"o kernel dictionary in the Hardy space (see \S 3 for a more general formulation).
 Using this $a_1$ in place of $a$ in (\ref{use1}), (\ref{use2}) and (\ref{use3}), we have that the corresponding reduced remainder $f_2$ has its least possible norm.
 To $f_2$ perform the same decomposition procedure, and so on. After $n$-iterations, we have
 \begin{eqnarray}
 \tilde{f}(z)=\sum_{k=1}^n\langle f_k,e_{a_k}\rangle B_k(z) + f_{n+1}(z)\prod_{k=1}^n\frac{z-a_k}{1-\overline{a}_kz},\end{eqnarray}
 where $\{B_k\}_{k=1}^\infty$ is the Takenaka-Malmquist system determined by
 $a_1,\cdots,a_k,\cdots,$ all in $\D,$ where
 \begin{eqnarray}\label{given} B_k(z)=e_{a_k}(z)\prod_{l=1}^{k-1}\frac{z-a_l}{1-\overline{a}_lz},\end{eqnarray}
 \begin{eqnarray}\label{msp} a_k=\max\{|\langle f_k,e_a\rangle|^2\ |\ a\in {\D}\},\end{eqnarray}
 \begin{eqnarray}\label{unavailable} f_{k+1}(z)=\frac{f_k(z)-\langle f_k,e_{a_k}\rangle e_{a_k}(z)}{\frac{z-a_k}{1-\overline{a}_kz}}\in H^2({\D}).\end{eqnarray}
We note that $\{B_k\}$ is automatically an orthonormal system, although not necessarily a basis. It turns out that under the maximal selections of $a_k, k=1,2,\cdots,$ there holds the convergence (\cite{QWa}):
\begin{eqnarray}
 \tilde{f}(z)=\sum_{k=1}^\infty\langle f_k,e_{a_k}\rangle B_k(z).\end{eqnarray}
 Due to the consecutive optimal selections of the parameters $a_k$ the convergence is in a fast pace. Although on the unit circle the Hardy space functions may not be smooth, it admits a promising convergence rate (\cite{QWa}).

  \begin{remark} Any sequence $(a_1,\cdots,a_n,\cdots)$ in $\D$ can define a TM system $\{B_k\}_{k=1}^\infty$ by (\ref{given}). A TM system is alternatively called a \emph{rational orthonormal system}. In the area of rational approximation, the study of TM systems together with their applications has a long history (\cite{Walsh}).  A TM system is an $H^p$-basis, $1<p<\infty,$ if and only if $\sum_{k=1}^\infty (1-|a_k|)=\infty.$ A half of the Fourier basis, $\{z^{k-1}\}_{k=1}^\infty,$ is a particular example of the basis cases. The study \cite{QWa} opens a new era of use of TM systems through adaptive selections of the parameters according to signals given in practical problems. The MSP of AFD declares the optimal selection principle at the one-step selection. It is the attainability of the
   global maximum at every one step that leads to necessity of repeating selections of some parameters, when it is the case. AFD shares the same idea as greedy algorithm for the one-step-optimal selection. The latter, however, does not address the issue of attainability of the global maximum in terms of parameters, nor addresses the necessity of repeating selections of the parameters, and nor concerns the multiple kernel concept (\cite{LT,MaZ,Te}). Through addressing those missed points in greedy algorithm AFD declares itself as a complete mathematical theory. AFD found close connections with the Beurling Theorem of the Hardy $H^2(\D)$ space asserting the direct-sum decomposition of the space into shift- and backward shift-invariant subspaces:
   \begin{eqnarray}\label{Beurling} H^2({\D})=\overline{\rm span}\{B_k\}_{k=1}^\infty \oplus \phi H^2({\D}),\end{eqnarray}
   where $\{B_k\}_{k=1}^\infty$ is the TM system defined by $a_1,\cdots,a_k,\cdots ,$ and $\phi$ is the Blaschke product, when can be defined. The sequence of complex numbers  can define a Blaschke product if and only if $\sum_{k=1}^\infty (1-|a_k|)<\infty.$ If the sequence cannot define a Blaschke product, then
   \begin{eqnarray}\label{cannot}H^2({\D})=\overline{\rm span}\{B_k\}_{k=1}^\infty.\end{eqnarray}
    With the AFD expansion we know that $\tilde{f}\in \overline{\rm span}\{B_k\}_{k=1}^\infty,$ the backward shift-invariant subspace in either the two cases (\ref{Beurling}) or (\ref{cannot}).
  \end{remark}

 \begin{remark} AFD was initially motivated by intrinsic positive phase derivative decomposition of analytic signals. It automatically generates a fast converging orthogonal expansion of which each entry has a meaningful instantaneous frequency. It has several variations, namely cyclic AFD, unwinding AFD, and has been generalized to multi-dimensions with the Clifford and several complex variables setting with scalar- to matrix-valued signals (\cite{Q2D,Qian2010,QWM,ACQS1,ACQS2}). In particular, the variation unwinding Blaschke expansion was first studied by Coifman, Steinerberger (\cite{CS}), and then Coifman and Peyri\'ere explored further connections with Blaschke products and outer functions (\cite{CP}). The type of decomposition was independently developed in \cite{Qian2010}. AFD has also been generalized to Hilbert spaces with a dictionary satisfying BVC (\cite{Q2D}, called \emph{pre-orthogonal adaptive Fourier decomposition} (POAFD). It, in particular, introduces the type of optimal sparse decomposition to general Hilbert spaces other than the Hardy type spaces (\cite{qu2018,qu2019}) in which there is no Blaschke-product-like functions. The Hilbert space generalizations are also for multivariate and non-scalar-valued functions. In a general Hilbert space, although there is no approximation theory as beautiful as one using TM system, POAFD can always produces a sparse series effectively approximating the given signal. We will give an exposition of POAFD in the latter half of the paper.  AFD and its variations, as well as its generalization POAFD, have become powerful tools in signal and image  analysis (\cite{Mi1,Mi2,CQT,LZQ,WQLG}).
\end{remark}

 \begin{remark}In the AFD algorithm, as a consequence of the orthogonality, there hold the relations:
\begin{eqnarray}\label{3}\langle f_k,e_{a_k}\rangle=\langle g_k,B_k\rangle=\langle \tilde{f},B_k\rangle,\ \quad k\ge 2, \end{eqnarray}
where
\begin{eqnarray}\label{SR} g_k(z)=\tilde{f}(z,w)-\sum_{l=1}^{k-1} \langle f_l,e_{a_l}\rangle B_l(z), \quad k\ge 2\end{eqnarray}
is the $k$-\emph{th} \emph{standard remainder}.
It is the relation (\ref{3}) that allows AFD to be generalized to Hilbert spaces with a dictionary satisfying BVC. In the latter there is no algebraic reduced remainder structure, nor explicit formulas for Gram-Schmidt orthogonalizations of parameterised reproducing kernels or dictionary elements that would play the role of TM systems being orthogonalizations of the Szeg\"o kernels in the classical Hardy  spaces (see the Appendix).\end{remark}

 Next we continue our theme on sparse representation of random signals. For an analytic random signal $f$ in $H^2_w({\D}),$ we obtain a sequence of parameters $a_1,a_2,\cdots,$ and an associated TM system $\{B_k\}_{k=1}^\infty$ that gives rise to an AFD sparse representation of the deterministic $\tilde{f}.$
 The question is: if we use the system $\{B_k\}_{k=1}^\infty$ to expand the original random signal
 $f(e^{it},w)=f_w(e^{it})$, then in what extent the expansion can represent the original RS $f?$ Or namely, what is the difference
 \begin{eqnarray}\label{d} d_f(e^{it},w)=f_w(e^{it})-\sum_{k=1}^\infty \langle f_w,B_k\rangle B_k(e^{it})?\end{eqnarray}
 We note that the RS $d_f$ has dependance on the TM system $\{B_k\}_{k=1}^\infty.$

In view of the Beurling Theorem, it well happens that for some $w$ the difference $d_f(e^{it},w)$ is non-zero. We have the following
\begin{theorem} \label{Sparse1}Let $f\in H^2_w({\D}),\tilde{f}=E_wf,$ and
$$\tilde{f}=\sum_{k=0}^\infty \langle \tilde{f},B_k\rangle B_k$$ be an AFD expansion of $\tilde{f}.$ Then, with the same $\{B_k\},$
\begin{eqnarray}\label{desired1}E_wd_f(e^{it},w)=0, \quad \forall t\in [0,2\pi).\end{eqnarray}
There holds the relation
\begin{eqnarray}\label{partial} E_w\|f_w-\sum_{k=1}^n\langle f_w,B_k\rangle B_k\|^2_{H^2_w}=\|d_f\|_{\N}^2+\sum_{k=n+1}^\infty
E_w|\langle f_w,B_k\rangle|^2,\end{eqnarray}
with \begin{eqnarray}\label{back1} \lim_{n\to\infty}\sum_{k=n+1}^\infty
E_w|\langle f_w,B_k\rangle|^2=0.\end{eqnarray}
\end{theorem}
And, in terms of the error $\r=f-\tilde{f}$ the difference, $d_f$ is estimated
\begin{eqnarray}\label{back2} \|d_f\|_{\N}^2=\|{r}\|_{\N}^2-\sum_{k=1}^\infty E_w|\langle {r}_w,B_k\rangle|^2.\end{eqnarray}

\noindent {\bf Proof} Since $\{B_k\}_{k=1}^\infty$ is an orthonormal system in the $\N$-space, the projection function $\sum_{k=1}^\infty \langle f_w,B_k\rangle B_k$ is in the Hilbert space $\N.$ The Bessel inequality gives
\[ \sum_{k=1}^\infty E_w|\langle f_w,B_k\rangle|^2\leq \|f\|_{\N}^2,\] that implies the desired relation (\ref{back1}).
As a consequence of the Riesz-Fisher Theorem the infinite series
\[\sum_{k=1}^\infty \langle f_w,B_k\rangle B_k\] is well defined
 for a. s. $w$ as a function in $H^2_w(\D).$ Hence the difference $d_f(w,\cdot)$  belongs to $H^2_w(\D).$ All these functions are in $\N.$

Since the underlying product measure space $\N$ is of finite total measure, both the convergence and the projection function are also in $L^1.$ As a consequence of the Fubini Theorem we can first take integral with respect to the probability, and get
\begin{eqnarray*}E_w(f_w-\sum_{k=1}^\infty \langle f_w,B_k\rangle B_k)&=&
\tilde{f}-E_w(\sum_{k=1}^\infty \langle f_w,B_k\rangle B_k)\\
&=&\tilde{f}-\sum_{k=1}^\infty E_w\langle f_w,B_k\rangle B_k\\
&=&\tilde{f}-\sum_{k=1}^\infty \langle \tilde{f},B_k\rangle B_k\\
&=&0,\end{eqnarray*}
as desired by (\ref{desired1}).

Noting that for each $w,$ $d_f$ is orthogonal with all $B_k$'s, we have the orthogonal decomposition
\[ f_w-\sum_{k=1}^n\langle f_w,B_k\rangle B_k=d_f+\sum_{k=n+1}^\infty
\langle f_w,B_k\rangle B_k,\]
that implies the desired Pythagoras relation (\ref{partial}).

Since
\[ d_f=(f_w-\tilde{f})-\sum_{k=1}^\infty \langle f_w-\tilde{f},B_k \rangle B_k={r}_w-\sum_{k=1}^\infty \langle {r}_w,B_k \rangle B_k,\]
estimating of $\|d_f\|^2_\N$ proceeds as
\begin{eqnarray*}
\|d_f\|^2_{\N}&=&E_w\int_0^{2\pi} |{r}_w(e^{it})-\sum_{k=1}^\infty \langle {r}_w,B_k\rangle B_k(e^{it})|^2dt\\
&=&E_w\left(\|{r}_w\|^2_{L^2}- \sum_{k=1}^\infty |\langle {r}_w,B_k\rangle|^2\right)\\
&=& \|{r}\|^2_{\N}-\sum_{k=1}^\infty E_w|\langle {r}_w,B_k\rangle|^2.\end{eqnarray*}
The proof of the theorem is complete.$\eop$\\

\begin{remark} The sparse random approximation corresponding to Theorem \ref{Sparse1} is designed for a deterministic signal with noise. Examples fitting into this theorem include those $r$ being white noise.
 In the following section we develop a sparse representation for analytic random signal that enjoys $d_f=0$ almost surely in $\Omega.$\end{remark}
 \bigskip

\subsection{The Type SAFDII: Taking Expectation Secondly}

\begin{theorem}\label{2MSP}
Let $f\in H^2_w({\D}).$ Then there exists $a_1\in {\D}$ such that
\begin{eqnarray}\label{msp} a_1=\arg \max \{E_w|\langle f_w,e_a\rangle|^2\ |\ a\in {\D}\}.\end{eqnarray}
\end{theorem}
\noindent{\bf Proof}\ \ Our effort will be rest on showing that the quantity under study satisfies a \emph{statistical boundary vanishing condition} (SBVC), that is
\begin{eqnarray}\label{BV}
\lim_{|a|\to 1}E_w|\langle f_w,e_a\rangle|^2=0.\end{eqnarray}
After proving (\ref{BV}) a density argument then concludes the desired maximal selection (\ref{msp}).
Since $f\in \N,$ the property
\begin{eqnarray}\label{indicate} E_w\sum_{k=0}^\infty |c_k(w)|^2<\infty\end{eqnarray}
implies that almost surely
\[ \sum_{k=0}^\infty |c_k(w)|^2<\infty.\]
As a consequence, almost surely $f_w(z)=\sum_{k=0}^\infty c_k(w)z^k \in H^2({\D}).$
Thanks to the BVC of the classical Hardy space (\cite{QWa}), we have almost surely
 \begin{eqnarray}\label{1}\lim_{|a|\to 1}|\langle f_w,e_a\rangle|^2=0.\end{eqnarray}
Now we show that there is a positive function of finite expectation dominating
 $|\langle f_w,e_a\rangle|^2$ a.s. in the process $|a|\to 1.$

In fact, for all $a\in \D$ uniformly
\begin{eqnarray*}
|\langle f_w,e_a\rangle|^2
\leq\|f_w\|^2
=\sum_{k=0}^\infty |c_k(w)|^2, \quad {\rm a.s}.
\end{eqnarray*}
The last positive random variable function, as a dominating function in $w,$ has a finite expectation as shown in (\ref{indicate}). Taking into account (\ref{1}), the Lebesgue domination convergence theorem can be used to conclude the desired SBVC (\ref{BV}). The proof is complete.$\eop$\\

The SAFDII algorithm proceeds as follows:
 Guaranteed by the Theorem \ref{2MSP}, with the same iterative steps as in AFD, one can select, at the $k$- step, an optimal $a_k:$
 \begin{eqnarray}\label{SMSP} a_k=\arg\max \{E_w|\langle (f_k)_w,e_a\rangle|^2\ |\ a\in {\D}\},\end{eqnarray}
 where $f=f_1,$ and
 \[ f_k(z,w)=(f_{k})_w(z)=\frac{(f_{k-1})_w(z)-\langle (f_{k-1})_w,e_{a_{k-1}}\rangle e_{a_{k-1}} (z)}{\frac{z-a_{k-1}}{1-\overline{a}_{k-1}z}}, \ \quad k\ge 2.\]
 The above maximal selection indicated in (\ref{SMSP}) is called \emph{stochastic maximal selection principle}, abbreviate as SMSP.
 We then construct a TM system $\{B_k\}_{k=1}^\infty,$ as given in (\ref{given}), corresponding to the selected $a_1,a_2,\cdots,$ and have the association
\[ f(z,w)\sim \sum_{k=1}^\infty \langle f_w,B_k\rangle B_k(z).\]
Due to the orthogonality of $\{B_k\}$ we also have
\begin{eqnarray}\label{OK}\langle (f_k)_w,e_{a_k}\rangle=\langle (g_k)_w,B_k\rangle=\langle f_w,B_k\rangle, \end{eqnarray}
where
\begin{eqnarray}\label{SR2} (g_k)_w(z)=g_k(z,w)=f(z,w)-\sum_{l=1}^{k-1} \langle f_w,B_l\rangle B_l(z), \quad k\ge 2, \end{eqnarray}
is the $k$-th standard remainder.
The relations (\ref{OK}) imply
\begin{eqnarray}\label{in} E_w|\langle (f_k)_w,e_{a_k}\rangle|^2=E_w|\langle (g_k)_w,B_k\rangle|^2=E_w|\langle f_w,B_k\rangle|^2.\end{eqnarray}
In view of (\ref{in}), the SMSP (\ref{SMSP}) is reduced to the form
\begin{eqnarray}\label{change form} a_k=\arg\max \{E_w|\langle f_w,B^a_k\rangle|^2\ |\ a\in {\D}\},\end{eqnarray}
 where
 \[ B^a_k(z)=e_a(z)\prod_{l=1}^{k-1}\frac{z-a_l}{1-\overline{a}_lz}.\]

We now prove
\begin{theorem}\label{N-convergence} Let $f(w,e^{it})\in H_w^2({\D})$ and $(a_1,\cdots,a_n,\cdots)$ be a sequence selected according to the SMSP given in (\ref{SMSP}), or, equivalently, (\ref{change form}). Then there holds, in the $\N$-norm sense,
 \begin{eqnarray}\label{convergence}f(z,w)=\sum_{k=1}^\infty \langle f_w,B_k\rangle B_k(z).\end{eqnarray}
 \end{theorem}
 \noindent{\bf Proof} By assuming the opposite we prove the convergence relation (\ref{convergence}) through a contradiction.  If the RHS does not converges to the LHS, then there is a non-trivial normal RS, $g\in \N,$ such that
 \begin{eqnarray}\label{contra}
 f(z,w)=\sum_{k=1}^\infty \langle f_w,B_k\rangle B_k(z)+g(z,w), \quad {\rm with}\ \quad \|g\|_{\N}> 0.\end{eqnarray}
 We note that $g$ is orthogonal with all $B_1,B_2,\cdots,B_k,\cdots,$ and
 \begin{eqnarray}\label{17}\|g\|_{\N}^2=\|f\|^2_{\N}-\sum^\infty_{k=1}E_w|\langle f_w,B_k\rangle|^2.\end{eqnarray}
 In particular,
 \begin{eqnarray}\label{part}
 \lim_{k\to \infty}E_w|\langle f_w,B_k\rangle|^2=0.
 \end{eqnarray}
 We show that there exists $b\in \D$ such that
 \begin{eqnarray}\label{positive} E_w|\langle g_w,e_b\rangle|^2>0.\end{eqnarray}
Denote $E_w|\langle g_w,e_b\rangle|^2=\delta^2.$
For, if (\ref{positive}) were not true, then almost surely for all $b\in \D$
\[ \langle g_w,e_b\rangle =0.\]
Due to the density of $\{e_b\}_{b\in {\D}}$ in $H^2({\D})$ we would have ,almost surely, $g_w=0$ as a function of $t,$ being contradictory to the condition $\|g\|_{\N}>0.$ We, in particular, can choose  $b$ being different from all the selected $a_k, k=1,2,\cdots$ We in below will fix this $b\in\D$ and proceed to derive a contradiction.

Set
\[ h_k=-\sum_{l=k}^\infty \langle f_w,B_l\rangle B_l.\]
From the definition of $g_k$ in (\ref{SR2}), there follows the orthogonal decomposition
\[ g=g_k+h_k.\]
The Bessel inequality implies, when $k$ is large,
\[ E_w|\langle h_k,e_b\rangle|^2\leq E_w\|h_k\|^2\leq \delta^2/4.\]
Hence
\[ 2E_w(|\langle g_k,e_b\rangle|^2)+\delta^2/2\ge E_w|\langle g_k,e_b\rangle +\langle h_k,e_b\rangle |^2=\delta^2,\]
which implies
\[ E_w|\langle g_k,e_b\rangle|^2\ge \delta^2/4.\]
Due to the reproducing kernel property of $e_b,$
for a large $k,$
\begin{eqnarray}\label{es}(1-|b|^2)^2E_w|(g_k)(b)|^2\ge \delta^2/4.\end{eqnarray}
Since pointwise there holds $f_k=g_k/\phi_k$ where
\begin{eqnarray}\label{depends} \phi_k(z)=\prod_{l=1}^{k}\frac{z-a_l}{1-\overline{a}_l}\quad {\rm and}\quad |\phi_k(b)|\leq 1, \quad \forall b\in {\D},\end{eqnarray}
there follows $|f_k|\ge |g_k|.$ Hence,
\[(1-|b|^2)^2E_w|(f_k)(b)|^2\ge \delta^2/4.\]
By using the reproducing property of $e_b$ again, the inner product form of the last equality has the form
\[ E_w|\langle f_k,e_b\rangle|^2=E_w|\langle f_w,B_k\rangle|^2\ge \delta^2/4\]
for all large enough $k.$  This is contradictory to (\ref{part}).
The proof is thus complete.$\eop$

\begin{remark} The original proof of convergence of AFD crucially depends on the property $\phi_k(b)\leq 1$ of Blaschke products (\cite{QWa}).   The succeeded generalizations of AFD on the contexts the matrix-valued Hardy space over the unit disc (\cite{ACQS1}), the Drury-Arveson space of several complex variables (\cite{ACQS2}), and the complex Hardy space over the $n$-torus (polydisc, \cite{Q2D}), are all based on
suitably defined Blaschke products with analogous unimodular properties like in (\ref{depends}).  As seen in the above proof, adaptability of the proof for the classical AFD to the stochastic version, Theorem \ref{N-convergence}, also rests on
(\ref{depends}). It can be further seen that the stochastic versions of the results in  \cite{ACQS1}, \cite{ACQS2} and \cite{Q2D} are available due to the same properties of their respective Blaschke products.
\end{remark}

\section{Stochastic SPOAFDs in Hilbert Spaces}

Our discussions on stochastic Hilbert spaces will refer to the existing study on deterministic Hilbert spaces with a dictionary satisfying BVC (see below). For  self-containing purpose we give a brief exposition on POAFD algorithm for deterministic signals (\cite{Q2D}, also see the survey paper \cite{CQT}).

\subsection{POAFD in a Hilbert Space With a Dictionary Satisfying BVC}
\def\C{\bf C}
The classical formulation of sparse representation of a Hilbert space is often under the assumption that the space has a dictionary that, by definition, is a dense subset of elements of the space of which each has unit norm. The unit norm requirement is not essential. We, therefore, release the norm-one requirement and only assume that the underlying Hilbert space, denoted by $\mathcal{H},$ has a dense subclass of elements $K_q, q\in \E,$ where the parameter set $\E$ is usually an open set of the complex plane, or  an open set of ${\R}^d$ and ${\C}^{d},$ or certain product spaces of them.  We denote the normalizations of $K_q$ by $E_q,$ i.e., $E_q=K_q/\|K_q\|, q\in {\E}.$ Below we often call the $K_q$'s by kernels. We now define what we call by $\lq\lq$multiple kernels". Let $(q_1,\cdots,q_n)$ be any $n$-tuple of parameters in $\E.$ Each of the terms $q_k, k=1,\cdots,n,$ may has multiplicity in the $k$-tuple $(q_1,\cdots,q_k).$ We denote by $l(k)$ the multiplicity of $q_k$ in $(q_1,\cdots,q_k).$  We accordingly introduce what we call \emph{multiple kernels} as follows. For simplicity we assume $\E\subset\C.$
For any $k\leq n,$ denote
\[\tilde{K}_{k}=\left[\left(\frac{\partial}{\partial \overline{q}}\right)^{(l(k)-1)}K_q\right](q_k),\]
where $l(k)$ is the multiple of $q_k$ in $(q_1,\cdots,q_k).$ With a little abuse of the notation, we will also denote $\tilde{K}_{k}$ by
 $\tilde{K}_{q_k}, k=1,2,\cdots,n,$ indicating the parameter sequence in use.  In complex-valued inner products due to use of conjugation of complex numbers the kernel $K_q$ is usually an anti-analytic function in $q$ and that is why the conjugate differential operator $\frac{\partial}{\partial \overline{q}}$ is used.
  The concept multiple kernel is a necessity
for \emph{pre-orthogonal maximal
 selection principle} (POMSP), stated as follows. Suppose we already have an $(n-1)$-tuple
 $\{q_1,\cdots,q_{n-1}\},$ with repetition or not, corresponding to
 the $(n-1)$-tuple $\{\tilde{K}_{q_1},\cdots,\tilde{K}_{q_{n-1}}\}.$
 By doing the G-S orthonormalization process consecutively we obtain
 an equivalent $(n-1)$-orthonormal basis
 $\{B_1,\cdots,B_{n-1}\}.$
 For any given $G$ in the Hilbert space we are interested in the supreme value
 \begin{eqnarray}\label{sup} \sup\{|\langle G,B_n^q\rangle| \ :\ q\in {\E}, q\ne q_1,\cdots,q_{n-1}\},\end{eqnarray}
where a finite supreme is guaranteed by the Cauchy-Schwartz inequality, and
$B_n^q$ is such that $\{B_1,\cdots,B_{n-1},B_n^q\}$ is the G-S orthonormalization
of $\{\tilde{K}_{q_1},\cdots,\tilde{K}_{q_{n-1}}, K_q\},$ where since $q\ne q_1,\cdots,q_{n-1},$ the function $B_n^q$ is precisely given by
 \begin{eqnarray}\label{GS}
 B_n^q
 =\frac{K_q-\sum_{k=1}^{n-1}\langle K_q,B_k\rangle_{\mathcal H}B_k}
 {\sqrt{\|K_q\|^2-\sum_{k=1}^{n-1}|\langle K_q,B_k\rangle_{\mathcal H}|^2}}.\end{eqnarray}
 There, however, may not exists $q_n\ne q_1,\cdots,q_{n-1}$ such that $B^{q_n}_n$ gives rise to the supreme of (\ref{sup}). The procedure in finding the supreme naturally leads to multiple kernels in relation to
what we call \emph{ Boundary Vanishing Condition} (\emph{BVC}):
 For any but fixed $G\in
\mathcal{H},$ if $p_n\in \E$ and $p_n\to\partial \E$ (including $\infty$ if
$\E$ is unbounded while in the case we use the compactification topology by adding the infinity point), then
\begin{eqnarray}\label{BVC} \lim_{n\to \infty} |\langle G,E_{p_n}\rangle|=0.\end{eqnarray}
Under the BVC given by (\ref{BVC}) a compact argument concludes that
 there exists a point $q_n\in \E$ and $q^{(l)}, l=1,2,\cdots,$
  such that $q^{(l)}$ are all different from $q_1,\cdots,q_{n-1},$
  $\lim_{l\to \infty}q^{(l)}=q_n,$ where $l(n)\ge 1,$ and
 \begin{eqnarray}\label{max} \lim_{l\to\infty}|\langle G,B_n^{q^{(l)}}\rangle|=
\sup\{|\langle G,B_n^q\rangle| \ :\ q\in {\E}, q\ne q_1,\cdots,q_{n-1}\}=|\langle G,B_n^{q_n}\rangle|,\end{eqnarray}
 where
 \begin{eqnarray}\label{GSmultiple}
 B_n^{q_n}
 =\frac{\tilde{K}_{q_n}-\sum_{k=1}^{n-1}\langle \tilde{K}_{q_n},B_k\rangle_{\mathcal H}B_k}
 {\sqrt{\|\tilde{K}_{q_n}\|^2-\sum_{k=1}^{n-1}|\langle \tilde{K}_{q_n},B_k\rangle_{\mathcal H}|^2}}.\end{eqnarray}
 Existence of such $q_n,$ as a consequence of BVC and multiple kernels, is the mathematical foundation of POAFD (\cite{Q2D,CQT}). We iteratively apply the above process to $G=G_n,$
 where $G_n$ is the standard remainder
\[ G_n=F-\sum_{k=1}^{n-1}\langle F,B_k\rangle B_k,\]
and $(B_1,\cdots,B_n)$ is the G-S orthogonalization of
$(\tilde{K}_{q_1},\cdots,\tilde{K}_{q_n}).$
Under the consecutive maximal selections of $\{q_k\}_{k=1}^\infty$ one eventually obtains, with a fast convergent pace,
\begin{eqnarray}\label{converge} F=\sum_{k=1}^\infty \langle F,B_k\rangle_{\mathcal{H}} B_k\end{eqnarray}
(\cite{Q2D,CQT}).

\begin{remark} We note that repeating selections of parameters can be avoided in practice. By definition of
supreme, for any $\rho\in (0,1),$ a parameter $q_n\in \E$
can be found, different from the previously
 selected $q_k, k=1,\cdots, n-1,$  to have
\begin{eqnarray}\label{rou} |\langle G_n,B_n^{q_n}\rangle|\ge \rho
\sup\{\langle G_n,B_n^q\rangle \ :\ q\in {\E}, q\ne q_1,\cdots,q_{n-1}\}.\end{eqnarray}
The corresponding algorithm for consecutively finding such a sequence  $\{q_n\}_{n=1}^\infty$  is called \emph{Weak Pre-orthogonal
Adaptive Fourier Decomposition} (WPOAFD). With WPOAFD one
may choose all $q_1,\cdots,q_n,\cdots
$ being distinguished. Under such selections we still get convergence (\ref{converge}) with a little less fast pace.\end{remark}
\begin{remark} An order $O(1/\sqrt{n})$  of the convergence rate can be proved: For $M>0,$ by defining

\begin{eqnarray}\label{Mcondition}
\mathcal{M}_M=\{ F\in \mathcal{H}: \exists \{c_n\}, \ \{E_{q_n}\} \
{\rm s.\ t.}\ F=\sum_{n=1}^\infty c_nE_{q_n}\ {\rm with}\  \sum_{n=1}^\infty |c_n|\leq M\},
\end{eqnarray}
for any $F\in \mathcal{M}_M,$ the POAFD partial sums satisfy
\[ \|F-\sum_{k=1}^n\langle F,B_k\rangle_{\mathcal H} B_k\|_{\mathcal H}\leq \frac{M}{\sqrt{n}}.\]
We note that the above convergence rate is the same as that of the Whittaker-Shannon interpolation formula in the sinc functions for bandlimited entire functions.
In the POAFD case the orthonormal system $\{B_1,\cdots,B_n,\cdots\}$ is not necessarily a basis but a system adapted to the given function $F.$ For the Hardy space case, due to the relations in (\ref{3}), the MSP (\ref{msp}) of AFD reduces to the MSP (\ref{max}) of POAFD, and thus AFD reduces to POAFD. The algorithm codes of AFD and POAFD, as well as those of some related ones are available at request (http://www.fst.umac.mo/en/staff/fsttq.html).\end{remark}

\begin{remark} AFD and POAFD have been seen to have two directions of developments. One is $n$-best kernel expansion. That is to determine $n$-parameters at one time, being obviously of better optimality in the $n$-sparse kernel approximation. The $n$-best approximation is motivated by the classical problem, yet still open in its ultimate global algorithm, called the best approximation to Hardy space functions by rational functions of degree not exceeding $n$ (\cite{Bara1986,Baratchart1991,QWM}).
The cyclic AFD and the gradient descending cyclic AFD (\cite{QWM}) may be adopted to give practical (not mathematical) $n$-best algorithms in Hilbert spaces with a dictionary satisfying BVC.
The second possible development of POAFD is related to exploration of Blaschke product-like functions and interpolation type problems in the general Hilbert context. For related publications see \cite{Qian2010,CS,ACQS1,ACQS2}. \end{remark}

\subsection{Stochastic POAFDs}
\def\H{\mathcal{H}}
Let $\mathcal{H}$ be a Hilbert space with a dense subset $\{K_q\}$ parameterized in an open set ${\E}: q\in \E.$ We assume that the dictionary satisfies BVC
\begin{eqnarray}\label{even more}
\lim_{q\to \partial \E} |\langle F,E_q\rangle|=0,
\end{eqnarray}
where $E_q=K_q/\|K_q\|.$  Let us consider random signals
$F(t,w), t\in T, w\in \Omega,$ where for a.s. $w\in \Omega,$ $F(\cdot,w)\in {\H};$ and for any $t\in T,$ $F(t,\cdot)$ is a random variable. Define
\begin{eqnarray}
\mathcal{N}(\mathcal{H},\Omega)&=&
\{F(t,w): F(\cdot,w)\in \mathcal{H},{\rm for\ a.s. \ w};  {\rm and}\ F(t,\cdot) \ {\rm being\ a\ random\ }\nonumber \\
& &\ \qquad{\rm variable\ for \ each\ fixed}\ t,\ {\rm and} \ E_w\|F(\cdot,w)\|^2_{\H}<\infty.\}
\end{eqnarray}
Denote $\|F\|_{{\N}({\H},\Omega)}=\left(E_w\|F(\cdot,w)\|^2_{\H}\right)^{1/2}.$ This formulation governs two types of stochastic POAFDs, abbreviated as SPOAFDI and SPOAFDII.

SPOAFDI is one to treat a noised deterministic signal by first taking the expectation and then doing maximal energy extractions.  We need to show $E_wF(t,w)\in {\H}.$  Following what is done in (\ref{similar}), by using the Minkovski inequality followed by the H\"older inequality, we get
 \[ \|E_wF(\cdot,w)\|_{\H}\leq E_w\|F(\cdot,w)\|_{\H}\leq \left(E_w\|F(\cdot,w)\|_{\H}^2\right)^{1/2}=\|F\|_{{\N}({\H},\Omega)}<\infty.\]
 This shows that the expectation belongs to the underlying Hilbert space $\H.$ Since $\H$ has a dictionary that satisfies BVC one can perform POAFD in $\H.$ The difference $d(t,w)=F(t,w)-E_wF(\cdot,w)$ enjoys the zero-expectation property and all the related quantities may be estimated as in the subsection 2.2. This approach gives rise to the type SPOAFDI that is suitable for analyzing signals corrupted with noise of zero expectation and of a small ${\N}(\mathcal{H},\Omega)$ norm.

 To perform the SPOAFDII type algorithm we first need to prove the stochastic boundary vanishing condition, or SBVC,
 \[\lim_{q\to\partial \E} E_w|\langle F_w,E_q\rangle|^2=0.\]
  To show this we, again, use the Lebesgue Dominated Convergence Theorem in the probability space, through showing \\

 1. First, for a.s. $w\in \Omega$
 \[ \lim_{q\to\partial \E}|\langle F_w,E_q\rangle|^2=0;\]
 and, \\
 2. secondly, for all $q$
 the function $|\langle F_w,E_q\rangle|^2$ is dominated by a positive integrable function in the probability space.\\
 \bigskip

 The property 1 is a consequence of BVC of the dictionary $\{E_q\}_{q\in \E}$ in $\H.$ To show 2, we have, by the H\"older inequality,
 \[ E_w|\langle F_w,E_q\rangle|^2\leq E_w\|F_w\|^2=\|F\|^2_{{\N}({\H},\Omega)}<\infty.\]
  This shows that $\|F_w\|^2$ is the dominating function in the probability space. The SBVC is hence proved. We have the following theorem

\begin{theorem}\label{SN-convergence} Let $F(t,w)\in {\N}({\H},\Omega)$ and $(q_1,\cdots,q_n,\cdots)$ be a consecutively selected parameter sequence according to SMSP
\[ q_k=\arg \max \{E_w|\langle {(G_k)}_w,B^q_k\rangle|^2\ |\ q\in {\E}\},\]
where
\[{(G_k)}_w=F_w-\sum_{l=1}^{k-1}\langle F_w,B_l\rangle B_l,\] and
$(B_1,\cdots,B_{k-1},B_k)$ is the G-S orthonormalization of $(B_1,\cdots,B_{k-1},\tilde{K}_{q_k}).$
Then there  holds, in the ${\N}({\H},\Omega)$-norm sense,
 \begin{eqnarray}\label{series}F(z,w)=\sum_{k=1}^\infty \langle F_w,B_k\rangle B_k(z).\end{eqnarray}
 \end{theorem}

 \begin{remark}The proof of Theorem \ref{N-convergence} crucially depends on the property
 $|\phi (z)|\leq 1$ of the classical Blaschke products. In the general Hilbert spaces case there may not exist Blaschke product-like functions playing the same role here, and, when there exist such functions, say $\phi,$ they may not enjoy the property $|\phi (z)|\leq 1.$ Below we give a proof of Theorem \ref{SN-convergence} that does not depend on Blaschke product-like functions. The proof is an adaptation of one for the deterministic signal case (see T. Qian, {\it A novel Fourier theory on non-linear phase and applications,} Advances in Mathematics (China), 2018, {47}(3), 321-347 (in Chinese), or the book chapter \cite{CQT}, the two references,  in two different languages, are essentially equivalent).  \end{remark}

 \noindent{\bf Proof of Theorem \ref{SN-convergence}}
  We will prove the theorem by contradiction.  If the RHS series of (\ref{series}) does not converges to the LHS function, then there is a non-trivial random signal $H\in {\N}({\H},\Omega)$ such that
 \begin{eqnarray}\label{contra}
 F(t,w)=\sum_{k=1}^\infty \langle F_w,B_k\rangle B_k(z)+H(z,w), \quad \|H\|_{{\N}({\H},\Omega)}> 0.\end{eqnarray}
 We note that $H$ is orthogonal with all $B_1,B_2,\cdots,B_k,\cdots,$ and
 \begin{eqnarray}\label{17}0<\|H\|_{{\N}({\H},\Omega)}^2=
 \|F\|^2_{{\N}({\H},\Omega)}-\sum^\infty_{k=1}E_w|\langle F_w,B_k\rangle|^2.\end{eqnarray}
 We claim that the fact $\|H\|_{{\N}({\H},\Omega)}>0$ implies that there exists $q\in \E,$ that determines some $\delta>0,$ such that
 \[ E_w|\langle H_w,E_q\rangle|^2=\delta^2>0.\]
For, if this were not true, then almost surely for all $q\in \E$
\[ \langle H_w,E_q\rangle =0.\]
Due to the density of $K_q$ in ${\N}({\H},\Omega)$ we would have almost surely $H_w=0$ as a function of $t,$ being contradictory to the condition $\|H\|_{{\N}({\H},\Omega)}>0.$ We, in particular, can choose  $q$ being distinguished from all the selected $q_k, k=1,2,\cdots.$ In below such $q\in \E$ will be fixed. The following argument will lead to a contradiction with the selections of $q_M$ for large enough $M.$

 Since $G_k$ is the $k$-standard remainders we rewrite the relation (\ref{contra}) as
\begin{eqnarray*} F_w&=&\left(\sum_{k=1}^M +\sum_{k=M+1}^\infty \right)\langle (G_k)_w,B_k\rangle B_k +H\\
&=&\sum_{k=1}^M\langle (G_k)_w,B_k\rangle B_k+\tilde{G}_{M+1}+H\\
&=&\sum_{k=1}^M\langle (G_k)_w,B_k\rangle B_k+G_{M+1},\end{eqnarray*}
where $$\tilde{G}_{M+1}=\sum_{k=M+1}^\infty \langle (G_k)_w,B_k\rangle B_k \quad {\rm and}\quad
G_{M+1}=\tilde{G}_{M+1}+H.$$
The Bessel inequality implies
\begin{eqnarray}\label{imply}
\lim_{M\to \infty}\|\tilde{G}_{M+1}\|_{{\N}({\H},\Omega)}=0.\end{eqnarray}
On one hand, we have, from (\ref{17}), for large enough $M,$
\begin{eqnarray}\label{Con} E_w|\langle (G_{M+1})_w,B_{M+1}\rangle|^2=E_w|\langle F_w,B_{M+1}\rangle|^2=E_w|\langle F_w,B^{q_{M+1}}_{M+1}\rangle|^2<\delta^2/16.\end{eqnarray}
On the other hand, we will show, for large $M,$ there holds
\begin{eqnarray}\label{hold}
E_w|\langle (G_{M+1})_w,B^q_{M+1}\rangle|^2>9\delta^2/16,\end{eqnarray}
where $B_{M+1}^q$
is the last function of the Gram-Schmidt orthonormalization of the $(M+1)$-system
$(B_1,B_2,\cdots,B_M,E_q)$ in the given order. The inequalities (\ref{Con}) and (\ref{hold}) the are contradictory.

So, all that remains is to show (\ref{hold}). From the triangle inequality of the ${\N}({\H},\Omega)$-norm,
\[ \left(E_w|\langle (G_{M+1})_w,B_{M+1}^q\rangle|^2\right)^{1/2}\ge
\left(E_w|\langle H_w,B_{M+1}^q\rangle|^2\right)^{1/2}-\left(E_w|\langle (\tilde{G}_{M+1})_w,B_{M+1}^q\rangle|^2\right)^{1/2}.\]
Using the Gauchy-Schwarz inequality and then (\ref{imply}), for large enough $M$ we have
\[ E_w|\langle (\tilde{G}_{M+1})_w,B_{M+1}^q\rangle|^2\leq \|\tilde{G}_{M+1}\|^2_{{\N}({\H},\Omega)}\leq \delta^2/16.\]
Therefore,
\begin{eqnarray}\label{combining} \left(E_w|\langle (G_{M+1})_w,B_{M+1}^q\rangle|^2\right)^{1/2}\ge
\left(E_w|\langle H_w,B_{M+1}^q\rangle|^2\right)^{1/2}-\delta/4.\end{eqnarray}
Next we compute the energy of the projection of $H_w$ into the span of $(B_1,\cdots,B_M,E_q).$
The energy is then
$E_w|\langle H_w,B_{M+1}^q\rangle|^2,$ for $H_w$ being orthogonal with
$B_1,\cdots,B_M.$ However,
the span is just the same if we alter the order $(B_1,\cdots,B_M,E_q)$ to
$(E_q,B_1,\cdots,B_M).$ As a consequence, the energy of the projection into the span is surely not less than the energy of $H_w$ projected onto the first function $E_q$ in the system. This gives rise to the relation
 \[E_w|\langle H_w,B_{M+1}^q\rangle|^2\ge E_w|\langle H_w,E_q\rangle|^2=\delta^2.\] Combining this result with (\ref{combining}), we have
\[ \left(E_w|\langle (G_{M+1})_w,B_{M+1}^q\rangle|^2\right)^{1/2}\ge 3\delta/4.\]
Thus we proved (\ref{hold}) that is contradictory with (\ref{Con}). This shows that the selection of $q_{M+1}$ did not obey SMSP, for we would better  select $q$ instead of $q_{M+1}$ at the $(M+1)$-th step. The proof of the theorem is hence
 complete.$\eop$

\begin{remark}
Both Theorem \ref{N-convergence} and Theorem \ref{SN-convergence} use orthonormal systems related to selected parameters according to the respective maximal selection principles. The difference is that in the AFD case one uses a kind of backward shift process, but in the POAFD case one uses the Gram-Schmidt orthonormalization. In Appendix we prove that, apart from a unimodular constant at each term, the TM system obtained from the backward shift process coincides with the orthonormal system obtained through the Gram-Schimidt process. Precisely, we will show
we will prove
\begin{theorem}\label{GS-TM}
Let $\{a_1,\cdots,a_n\}$ be any $n$-tuple of parameters
in $\D$ in which multiplicities are allowed.
Denote by $l(m)$ the multiplicity of $a_m$ in the $m$-tuple $\{a_1,\cdots,a_{m}\}, 1\leq m\leq n.$ For each $m,$ denote by \[\tilde{k}_{a_m}(z)=\frac{\partial^{l(m)-1}}{(\partial \overline{a})^{l(m)-1}}k_a(z)|_{a=a_m},\
{\rm where}\  k_a(z)=\frac{1}{1-\overline{a}z}.\]
Then the Gram-Schmidt orthonormalization of $\{\tilde{k}_{a_1},\cdots,\tilde{k}_{a_m}\}$ in the given order coincides with the $m$-TM system $\{B_1,\cdots,B_m\}$
(\ref{given}) defined through the ordered $m$-tuple $\{a_1,\cdots,a_m\}.$
\end{theorem}
Based on this result Theorem \ref{N-convergence} is, as a matter of fact, a corollary of Theorem \ref{SN-convergence}.
\end{remark}

\section{Acknowledgement}
The author wishes to express his sincere thankfulness to Dr Chen Wei-Guo, Dr Wang Shi-Lin, Dr Cheng Han-Sheng, Prof Chen Qiu-Hui, Prof Leong Ieng Tak for their interest and encouragement to study this topic, and useful comments on draft material of this article.

\section{Appendix}

\noindent{\bf Proof of Theorem \ref{GS-TM}}
Denote the canonical Blaschke product determined by $a_1,\cdots,a_m$ as
\[\phi_{a_1,\cdots,a_m}(z)=\prod_{l=1}^m\frac{z-a_l}{1-\overline{a}_lz}.\]
We first show that for any $a\in \D$ being different from $a_1,\cdots,a_{m-1}$ there holds
\begin{eqnarray}\label{induc}
k_a(z)-\sum_{l=1}^{m-1}\langle k_a,B_l\rangle B_l(z)=\overline{\phi}_{a_1,\cdots,a_{m-1}}(a)\phi_{a_1,\cdots,a_{m-1}}(z)k_a(z).
\end{eqnarray}
For this aim we use mathematical induction. First we verify the case $m=2.$ Using the reproducing kernel property of $k_a,$ there follows
\begin{eqnarray*}
k_a-\langle k_a,B_1\rangle B_1(z)&=&\frac{1}{1-\overline{a}z}-\overline{B}_1(a)B_1(z)\\
&=&\frac{1}{1-\overline{a}z}-\frac{\alpha}{1-\overline{a}_1z},\qquad \alpha=\frac{1-|a_1|^2}{1-a_1\overline{a}},\\
&=& \frac{\overline{a}-a_1}{1-a_1\overline{a}}\frac{z-a_1}{1-\overline{a}_1z}
\frac{1}{1-\overline{a}z}\\
&=&\overline{\phi}_{a_1}(a)\phi_{a_1}(z)k_a(z).\end{eqnarray*}
 Assume that (\ref{induc}) holds for $m$ being replaced by $m-1.$ Under this inductive hypothesis, we have
 \begin{eqnarray*}
 k_a(z)-\sum_{l=1}^{m-1}\langle k_a,B_l\rangle B_l(z)&=&[ k_a(z)-\sum_{l=1}^{m-2}\langle k_a,B_l\rangle B_l(z)]-\langle k_a,B_{m-1}\rangle B_{m-1}(z)\\
 &=&\overline{\phi}_{a_1,\cdots,a_{m-2}}(a)\phi_{a_1,\cdots,a_{m-2}}(z)k_a(z)-\langle k_a,B_{m-1}\rangle B_{m-1}(z)\\
 &=&\overline{\phi}_{a_1,\cdots,a_{m-2}}(a)\phi_{a_1,\cdots,a_{m-2}}(z)k_a(z)-
\overline{B}_{m-1}(a)B_{m-1}(z)\\
&=&\overline{\phi}_{a_1,\cdots,a_{m-2}}(a)\phi_{a_1,\cdots,a_{m-2}}(z)k_a(z)
\left[k_a(z)-\frac{1-|a_{m-1}|^2}{(1-a_{m-1}\overline{a})(1-\overline{a}_{m-1}z)}\right]\\
&=&\overline{\phi}_{a_1,\cdots,a_{m-1}}(a)\phi_{a_1,\cdots,a_{m-1}}(z)k_a(z).
\end{eqnarray*}
We hence proved (\ref{induc}). Next we deal with the orthonormalization allowing repetition of the parameters. Now we are with the new inductive hypothesis that the Gram-Schmidt orthonormalization of $\{\tilde{k}_{a_1},\cdots,\tilde{k}_{a_{m-1}}\}$ is the $(m-1)$-TM system $\{B_1,\cdots,B_{m-1}\}.$ First assume $a_m$ is different from all the preceding $a_k, k=1,\cdots,m-1.$ In (\ref{induc}) let $a=a_m.$ By taking the norm on the both sides of (\ref{induc}) and invoking the orthonormality of the TM system we have
\[ \|k_{a_m}(z)-\sum_{l=1}^{m-1}\langle k_{a_m},B_l\rangle B_l(z)\|=e^{-ic}\overline{\phi}_{a_1,\cdots,a_{m-1}}(a_m)
\frac{1}{\sqrt{1-|a_m|^2}},\]
where $c$ is a real number depending on $a_m$ and $a_1,\cdots,a_{m-1}.$
We thus conclude that
\begin{eqnarray}\label{Bla} \frac{k_{a_m}(z)-\sum_{l=1}^{m-1}\langle k_{a_m},B_l\rangle B_l(z)}{\|k_{a_m}(z)-\sum_{l=1}^{m-1}\langle k_{a_m},B_l\rangle B_l(z)\|}=e^{ic}\phi_{a_1,\cdots,a_{m-1}}(z)e_{a_m}(z).\end{eqnarray}
Note that here we have the case $k_{a_m}=\tilde{k}_{a_m}$ and $l(m)=1.$ Next we extend the above relation to the cases that $a=a_m$ coincides with some of the preceding $a_1,\cdots,a_{m-1}.$ In that case we have $l(m)>1,$ and we are to show
\begin{eqnarray}\label{tildeBla} \frac{\tilde{k}_{a_m}(z)-\sum_{l=1}^{m-1}\langle \tilde{k}_{a_m},B_l\rangle B_l(z)}{\|\tilde{k}_{a_m}(z)-\sum_{l=1}^{m-1}\langle \tilde{k}_{a_m},B_l\rangle B_l(z)\|}=e^{ic}\phi_{a_1,\cdots,a_{m-1}}(z)e_{a_m}(z),\end{eqnarray}
where $c$ depends on $a_1,\cdots, a_m.$
For $b$ being sufficiently close to $a_m$ in $\D$ we have up to the
$(l(m)-1)$-order power series expansion in the variable $b:$
\begin{eqnarray*} k_b(z)
&=&\sum_{l=0}^{l(m)-1}\frac{1}{l!}
\left[\frac{\partial}{\partial\overline{a}}\right]^{l}k_a(z)|_{a=a_m}(b-a_m)^{l}+
o((b-a_m)^{(l(m)-1)}\\
&=&T(z)+\frac{1}{(l(m)-1)!}
\tilde{k}_{a_m}(z)(b-a_m)^{l(m)-1}+
o((b-a_m)^{(l(m)-1)},\end{eqnarray*}
where
\[T(z)=\sum_{l=0}^{l(m)-2}\frac{1}{l!}
\left[\frac{\partial}{\partial\overline{a}}\right]^{l}k_a(z)|_{a=a_m}(b-a_m)^{l}.\]
Now, according to the inductive hypothesis, $B_1,\cdots,B_{m-1}$ involve the derivatives of the reproducing kernel up to the $(l(m)-2)$-order, and hence
\begin{eqnarray}\label{insert} T(z)-\sum_{k=1}^{m-1}\langle T,B_k\rangle B_k=0.\end{eqnarray}
Inserting the left-hand-side of (\ref{insert}) into (\ref{Bla}), where $a_m$ is replaced by $b$ with $b\to a_m$ horizontally (meaning that ${\rm Im}(b)={\rm Im}(a_m)$), while dividing by $(b-a_m)^{l(m)-1}>0,$ we have
\begin{eqnarray*}
\frac{\frac{k_b(z)-T(z)}{(b-a_m)^{l(m)-1}}-\sum_{l=1}^{m-1}\langle \frac{k_b-T}{(b-a_m)^{l(m)-1}},B_l\rangle B_l(z)}{\|\frac{k_b-T}{(b-a_m)^{l(m)-1}}-\sum_{l=1}^{m-1}\langle \frac{k_b-T}{(b-a_m)^{l(m)-1}},B_l\rangle B_l(z)\|}=e^{ic_b}
\phi_{a_1,\cdots,a_{m-1}}(z)k_b(z).
\end{eqnarray*}
Letting $(b-a_m)\downarrow 0$ and noticing that the Taylor series remainder is an infinitesimal of an order higher than $(b-a_m)^{l(m)-1}$,
 we obtain the desired relation (\ref{tildeBla}).
The proof is complete.$\eop$\\
The above is a constructive proof. Another proof may be found in Qian and Wegert 2013 (also see B. Ninness et al. \cite{NG}, 1997).
As far as what the author is aware of, the unit disc and a half of the complex plane are the only cases to which the equivalence of the two processes, i.e., the Blaschke product-backward shift formulation and the Gram-Schmidt orthogonalization, has been proved.

\end{document}